\documentclass[psamsfonts]{conm-p-l}
\usepackage{amsmath,amssymb,euscript,survey}

\newtheorem{theorem}{Theorem}[section]
\newtheorem{lemma}[theorem]{Lemma}
\newtheorem{proposition}[theorem]{Proposition}
\newtheorem{corollary}[theorem]{Corollary}

\theoremstyle{definition}
\newtheorem{definition}[theorem]{Definition}

\theoremstyle{remark}
\newtheorem{remark}[theorem]{Remark}

\numberwithin{equation}{section}

\begin{document}

\title{Approximation in ergodic theory, Borel, and Cantor dynamics}

%    Information for first author
\author{S.~Bezuglyi}
%    Address of record for the research reported here
\address{Department of Mathematics, Institute for Low Temperature Physics,
47 Lenin ave., 61103 Kharkov, Ukraine}
%    Current address
\email{bezuglyi@ilt.kharkov.ua}
%    \thanks will become a 1st page footnote.
\thanks{S. Bezuglyi was supported in part by CRDF grant UM1-2546-KH-03}

%    Information for second author
\author{J.~Kwiatkowski}
\address{Faculty of Mathematics and Computer Sciences, Nicolas Copernicus
University, ul. Chopina 12/18, 87 - 100 Torun, Poland}
\email{jkwiat@mat.uni.torun.pl}
\thanks{J. Kwiatkowski was supported in part by KBN grant 1P 03A 038 26}

% Third author

\author{K. Medynets}
\address{Department of Mathematics, Institute for Low Temperature Physics,
47 Lenin ave., 61103 Kharkov, Ukraine}
\email{medynets@ilt.kharkov.ua}

%    General info
\subjclass{Primary 37A40, 37B05, 03066}
\date{}

%\dedicatory{This paper is dedicated to our advisors.}

\keywords{automorphism, measure space, standard Borel space,
Cantor set, homeomorphism}

\newcommand{\C}{\mathcal C}
\newcommand{\B}{\mathcal B}
\newcommand{\aut}{Aut(X,\B)}
\newcommand{\per}{{\mathcal P}er}
\newcommand{\ap}{{\mathcal A}p}
\newcommand{\meas}{{\mathcal M}_1(X)}
\newcommand{\Aut}{Aut(X,\B,\mu)}
\newcommand{\N}{\mathbb N}
\newcommand{\Z}{\mathbb Z}
\newcommand{\e}{\varepsilon}
\newcommand{\bs}{(X,\B)}
\newcommand{\h}{Homeo(\Omega)}
\newcommand{\Om}{\Omega}

\begin{abstract}
This survey is focused on the results related to topologies on the groups
of transformations in ergodic theory, Borel, and Cantor dynamics. Various
topological properties (density, connectedness, genericity) of these groups
and their subsets (subgroups) are studied.
\end{abstract}

\maketitle

In this paper, we intend to present a unified approach to the study of
topological properties of transformation groups arisen in ergodic theory,
Borel, and Cantor dynamics. We consider here the group $\Aut$ (and
$Aut_0(X,\B,\mu))$ of all non-singular (measure preserving, resp.)
automorphisms of a standard measure space, the group $\aut$ of all Borel
automorphisms of a standard Borel space, and the group $\h$ of all
homeomorphisms of a Cantor set $\Omega$. The basic technique in the study
of transformation groups acting on an underlying space is to introduce
various topologies into these groups which make them topological groups and
investigate topological properties of the groups and their subsets
(subgroups).

The study of topologies on the group of transformations of a space has a
long history. The first significant results in this area are the classical
results of Oxtoby and Ulam on the typical dynamical behavior of
homeomorphisms which preserve a measure \cite{O-U}. This circle of problems
has attracted attention in various areas of dynamical systems, notably, in
measurable and topological dynamics, where it is important for many
applications to understand what kind of transformations is typical for
certain dynamics. Of course, this problem assumes that a topology is
defined on the group of all transformations.

The topologies which are usually considered on the groups $\Aut$ and
$Aut_0(X,\B,\mu)$ were defined 60 years ago in the pioneering paper by
Halmos \cite{Hal 1} where he called them the uniform and weak topologies.
The use of these topologies turned out to be very fruitful and led to many
outstanding results in ergodic theory. For instance, the best known
theorems concerning ergodic, mixing, and weakly mixing automorphisms of a
measure space were obtained by P.~Halmos and V.A.~Rokhlin (see \cite{Hal
1}, \cite{Hal 2}, \cite{Ro 1}, \cite{Ro 2}). Many further statements on
approximation of automorphisms of a measure space can be found in the book
\cite{C-F-S} and in a great number of research and expository papers
written in last decades (see, e.g., the references to this survey).

Motivated by ideas used in ergodic theory, we define several topologies on
the groups $\aut$ and $\h$ which are similar to the weak and uniform
topologies. They were first introduced in the context of Cantor dynamics in
\cite{B-K 1} and \cite{B-K 2} where we marked that many results in  Cantor
dynamics have their counterparts ergodic theory.  Then, we used these
topologies to the study of the groups $\aut$ and $\h$ in the series of
papers \cite{B-D-K 1}, \cite{B-D-K 2}, \cite{B-D-M}, \cite{B-M}. To the
best our knowledge, topologies on the group $\aut$ had not been considered
before. On the other hand, there exist a few useful topologies on the group
of homeomorphisms of a topological space $K$. The most known of them is the
topology of uniform convergence on $Homeo(K)$ (or the compact-open
topology). Many remarkable results about topological properties of
$Homeo(K)$ were obtained for  a connected compact metric space $K$ (see,
e.g. \cite{Alp-Pr 2} and the references there).

In our study of the groups $\aut$ and $\h$,  we noticed that most of
topological properties known for $\Aut$ hold also for  Borel and Cantor
dynamics. Based on this observation, we have tried to systematize the known
results for measurable,  Borel, and Cantor dynamics to emphasize some
common features and point out the existing distinguishes. In the article,
we consider the global topological structure of these large groups. We also
study some classes (subsets, subgroups) of transformations which are
naturally defined, for instance, periodic, aperiodic, minimal, odometers,
etc. The closures of these classes are found in the topologies we have
defined. We consider here the problems  and results which are parallel in
ergodic theory, Borel, and Cantor dynamics. It will be seen below that the
groups we study have a lot of similar topological properties. This
observation is a justification for studying topologies on the groups of
transformations in an unified manner.

The paper is organized as follows. In the first section, we collect main
definitions and notions used in ergodic theory, Borel, and Cantor dynamics.
Our goal is two-fold: we fix the notation  used in the paper and try to
make the paper self-contained (this may be useful for non-specialists).
Section 2 contains the well known results on applications of the weak and
uniform topologies in ergodic theory. We do not intend to give a
comprehensive description of all results related to topological properties
of $\Aut$ and $Aut_0(X,\B,\mu)$. They are too numerous to list here. Our
choice is explained  by interest and knowledge, and, all above, by our
study of topologies in the context of Borel and Cantor dynamics. The next
two sections contain the results on topologies on the groups $\aut$ and
$\h$. We define two topologies $\tau$ and $p$ which are analogous to the
uniform and weak topologies in measurable  dynamics and formulate many
results on topological properties of these groups and various subsets and
subgroups. In particular, we consider periodic, aperiodic, minimal
transformations, odometers etc. In the last section, we give a comparison
of results obtained in measurable,  Borel, and Cantor dynamics.

While preparing this survey, we  improved several results from our previous
works and included them in the paper with sketches of the proofs. We also
used the nice written papers \cite{Ch-Pr}, \cite{Alp-Pr 2} where one can
find many results on topological properties of the group of automorphisms
of a measure space and the group of measure-preserving homeomorphisms of a
compact connected manifold.

In this survey, we do not discuss many remarkable results which have been
published during last decade in Borel and Cantor dynamics. In this
connection, we refer to \cite{G-P-S 1}, \cite{G-P-S 3}, \cite{Hj},
\cite{J-Ke-L}, \cite{Ke-Mil} where further references can be found.

\section{Preliminaries}

In this section, we recall main notions and definitions from ergodic theory
(called also measurable dynamics), Borel, and Cantor dynamics which are
used in the paper. We will keep the notation from this section throughout
the paper.

%%%%%%%%%%% Borel dynamics %%%%%%%%%%%%%%%%%%%%%%%%

\subsection{Borel dynamics.}\label{Bor-Dyn} Let $X$ be a topological space.
The class $\mathcal B(X) = \B$ of Borel sets of $X$ is the $\sigma$-algebra
generated by open sets. We call $(X,\B)$ the Borel space. A separable
complete metric space $Y$ is called a Polish space. Let ${\mathcal C}$ be
the $\sigma$-algebra of Borel subsets of $Y$. A Borel space $(X,\B)$ is
called a {\it standard Borel space} if it is isomorphic to a Polish space
$(Y,\C)$. Any two standard Borel spaces are isomorphic.

Let $\meas$ denote the set of all Borel probability measures on $(X,\B)$.
A measure $\mu\in \meas$ is called {\it non-atomic} (or continuous) if
$\mu(\{x\}) = 0\ \forall x\in X$. The measure, supported by a point $x\in
X$, is denoted by $\delta_x$. Two measures $\mu$ and $\nu$ are called
{\it equivalent} (in symbols, $\mu \sim \nu$) if they have the same
collection of subsets of measure zero.

A one-to-one Borel map $T$ of $X$ onto itself is called an {\it
automorphism} of $(X,\B)$. Denote by $\aut$ the group of all Borel
automorphisms of $(X,\B)$. For $T\in Aut(X,B)$, the set $supp(T)=\{x\in X
: Tx\neq x\}$ is called the {\it support} of $T$. Let $Ctbl(X)$ denote the
set of automorphisms with at most countable support, i.e., $T\in Ctbl(X)$
if $|(supp (T))|\leq \aleph_0$. Clearly, $Ctbl(X)$ is a normal subgroup of
$Aut(X,B)$. The quotient group $\widehat{Aut}(X,B)= Aut(X,B)/Ctbl(X)$ is
obtained  by identifying automorphisms from $\aut$ which differ on an at most
countable set.

For $T\in \aut$ and $x\in X$, let $Orb_T(x)=\{T^nx: n\in\mathbb Z\}$ denote
the $T$-orbit of $x$. A point $x\in X$ is called a {\it periodic point} of
period $n$ if  $T^nx=x$ and $T^ix \neq x,\ i=1,...,n-1$. An automorphism
$T$ is called  {\it pointwise periodic} (or periodic) if every point $x\in
X$ is periodic for $T$, and $T$ is called {\it aperiodic} if it has no
periodic points, i.e. every $T$-orbit is infinite. We denote by $\per$ and
$\ap$ the sets of pointwise periodic and aperiodic automorphisms,
respectively.

A Borel automorphism $T\in Aut(X,\B)$ is called {\it smooth} if there is a
Borel subset $A\subset X$ (called {\it wandering}) that meets every
$T$-orbit exactly ones. We denote the class of smooth automorphisms by
$\mathcal{S}m$. By our definition, every pointwise periodic automorphism is
smooth.

Given $T\in \aut$, denote by $M_1(T)$ the set of $T$-invariant Borel
probability measures. In particular, $M_1(T)$ may be empty. An automorphism
$T$ is called {\it incompressible} if $M_1(T) \neq \emptyset$ (in symbols,
$T\in \mathcal{I}nc$).

By definition, the {\it full group} of $T\in Aut(X,\B)$ is
$[T]=\{S\in Aut(X,\B) : Orb_S(x)\subset Orb_T(x),\;x\in X\}$.

We recall also the definition of odometers.  Let $\{p_n\}_{n=0}^\infty$ be
a sequence of integers such that $p_n\geq 2$. Let
$\Omega=\prod_n\{0,1,\ldots,p_n-1\}$ be endowed with the product topology.
Then $\Omega$ is a Cantor set. Let $S : \Omega\rightarrow \Omega$ be
defined as follows: $S(p_0-1,p_1-1,\ldots)=(0,0,\ldots)$, and for any other
$x= (x_k)\in X$, find the least $k$ such that $x_k\neq p_k-1$ and set
$S(x)=(0,0,\ldots,0,x_k+1,x_{k+1},x_{k+2},\ldots)$. More general, a Borel
automorphism $T$ is called an  {\it odometer} if it is Borel
isomorphic to some $S$. In other words, there exists a sequence of
refining partitions of $X$ such that each partition is a $T$-tower and the
elements of the partitions generate $\B$. We will denote by $\mathcal Od$
the set of odometers.

%%%%%%%%%%%%%%%%%%%Measurable dynamics%%%%%%%%%%%%

\subsection{Measurable dynamics.}\label{Meas-Dyn}
A {\it measure space} $(X, \mathcal A, \mu)$ consists of a set $X$, a
$\sigma$-algebra $\mathcal A$ of subsets of $X$, and a measure (finite or
$\sigma$-finite) $\mu$ on $\mathcal A$ with respect to which $\mathcal A$
is complete. If $T$ is a one-to-one measurable map of $X$ onto itself, then
we can define the measure $\mu\circ T(A) := \mu(TA),\; A\in \mathcal A$. If
$\mu\circ T \sim \mu$, then $T$ is called a {\it non-singular} automorphism
of $(X, \mathcal A, \mu)$. If $\mu\circ T = \mu$, then $T$ is called {\it
measure-preserving}. The group of all non-singular automorphisms is denoted
by $Aut(X, \mathcal A, \mu)$, and that of measure-preserving automorphisms
is denoted by $Aut_0(X, \mathcal A, \mu)$. As usual in ergodic theory, we
use the ``mod 0 convention'', that is all equalities between sets,
automorphisms, etc. are to be understood modulo null sets.

We say that $(X,\B,\mu)$ is a {\it standard measure space} (or {\it
Lebesgue space}) if it is isomorphic to $([0, a], \lambda)$ where $a\in
[0,\infty]$ and $\lambda$ is the Lebesgue measure.\footnote{In contrast to
the original definition, we consider here only non-atomic measures on
$X$.} If $\mu(X) =1$, then $(X,\B,\mu)$ is called a probability standard
measure space. To illustrate these definitions, we note that any
uncountable Borel subset of a Polish space together with a Borel  measure
on it is an example of a standard measure space.

An automorphism $T\in \Aut$ is called {\it ergodic} (or $\mu$-ergodic) if
for every $T$-invariant set $A\in \B$ one has either $\mu(A) = 0$ or $\mu(X
\setminus A)=0$. An automorphism $T\in \Aut$ is called {\it mixing} if
$\mu(T^{-k}A\cap B)\to \mu(A)\mu(B)$ as $k\to\infty$  for all $A,B\in\B$.
$T$ is {\it weakly mixing} if $T \times T$ is $\mu\times \mu$-ergodic on
$X\times X$. The {\it full group} $[T]$ generated by $T\in \Aut$ consists
of those automorphisms whose orbits are contained in $T$-orbits. By
definition, the {\it normalizer} $N[T]$ of $[T]$ is the set $\{S\in \Aut :
S[T]S^{-1} = [T]\}$, and the centralizer $C(T)$ of $T$ is the set $\{S\in
\Aut : ST = TS\}$.

Recall that an automorphism $T\in \Aut$ is called of type $II$ if
there exists a $T$-invariant measure $\nu$ equivalent to $\mu$. If
$\nu$ is finite, then $T$ is of type $II_1$; if $\nu$ is infinite
then $T$ is of type $II_\infty$. If $T$ does not admit an invariant
(finite or infinite) measure equivalent to $\mu$, then $T$ is
called of type $III$. Automorphisms of type $III$ are divided into
the subtypes $III_0$, $III_\lambda\ (0< \lambda < 1)$, and $III_1$
(see \cite{Ha-Os} for details).

%%%%%%%%%%%%%%%%%%%%%Topological dynamics%%%%%%%%%%%%%%%%%%%

\subsection{\bf Topological dynamics.}\label{Top_Dyn}
Let $K$ be a compact metric space. Denote by $Homeo(K)$ the group of all
homeomorphisms of $K$. We say that $T\in Homeo(K)$ is {\it topologically
transitive} if there is $x\in K$ such that $Orb_T(x)$ is dense in $K$.
Also, $T$ is called {\it minimal} if every $T$-orbit is dense in $K$ or,
equivalently, $T$ has no proper closed invariant subsets. It is  said that
$T\in Homeo(K)$ is {\it mixing} if for every non-empty open subsets $U,V$
there exists $N\in \N$ such that $T^nU \cap V \neq\emptyset$ for all
$n\geq N$. The sets of all minimal and mixing homeomorphisms will be
denoted by $\mathcal Min$ and $\mathcal Mix$, respectively. A
homeomorphism $T$ of $K$ is called {\it weakly mixing} if $T\times T$ is
transitive on $K\times K$. Notice that an odometer $S$ is obviously a
minimal homeomorphism.

We say that $\Omega$ is a {\it Cantor set} if it is homeomorphic to a
compact metric space without isolated points whose base consists
of clopen sets. All Cantor sets are homeomorphic.

For aperiodic $T\in Homeo(\Omega)$, define the {\it full group of $T$} as
$[T]=\{S\in Homeo(\Omega) : Orb_S(x)\subset Orb_T(x)\ \mbox{ for all}\ x\in
\Omega\}$. Thus, if $S\in [T]$, then $Sx=T^{n_S(x)}x$ for all $x$, and
$x\mapsto n_S(x)$ is a Borel map from $\Omega$ to $\mathbb Z$. The
countable subgroup $[[T]]=\{S\in[T] : n_S(x)\mbox{ is continuous}\}$ is
called the {\it topological full group} of $T$.

%%%%%%%%%%%%%%%%%%%%%Section 2%%%%%%%%%%%%%%%%%%%%%%%%

\section{The uniform and weak topologies on $Aut(X,\B,\mu)$}

In this section, we collect a number of results about topological
properties of $\Aut$ and $Aut_0(X, \B,\mu)$. These results are not new and
presented here to emphasize some parallels with results proved for Borel
and Cantor dynamics.

\subsection{\bf Topological properties of the space $\Aut$.}
Let $Aut(X,\B,\mu)$ be the group of all non-singular measurable
automorphisms of a standard measure space $(X,\B,\mu)$ with
$\mu(X) = 1$. Using the ideas which first appeared in \cite{Hal
1}, one can introduce two topologies on $Aut(X,\B,\mu)$ generated
by the metrics $d_u$ and  $d_w$. For $S,T \in Aut(X,\B,\mu)$, set
\begin{equation}\label{d_u}
    d_u(S,T) = \mu(\{x\in X : Sx\neq Tx\}),
\end{equation}
\begin{equation}\label{d_w}
    d_w(S,T) = \sum_{n=1}^\infty 2^{-n}\mu(SA_n \Delta TA_n),
\end{equation}
where $(A_n)$ is a countable family of measurable subsets generating $\B$.

It is well known that the topologies defined by the metrics (\ref{d_u})
and (\ref{d_w}) do not depend on a measure from the class $[\mu]$ of
measures equivalent to $\mu$. We will follow the long established tradition
in ergodic theory to call them the {\it uniform and weak topologies},
respectively. Clearly, the uniform topology is strictly stronger than the
weak topology. These topologies on the groups $Aut(X,\B,\mu)$ and
$Aut_0(X,\B,\mu)$ have been studied in many papers. We summarize the obtained
results on topological and algebraic properties of these groups in the
next two theorems.

\begin{theorem}[uniform topology]\label{Aut-un}
$(1)$ The uniform topology on $Aut(X,\B,\mu)$ is equivalent to the
topology generated by the  metric\footnote{This fact is a justification of
the name  `uniform' for this topology.}
\begin{equation}\label{d'_u}
    d'_u(S,T) = \sup_{A\in \B} \mu(SA \Delta TA).
\end{equation}
$(2)$ $(Aut(X,\B,\mu), d_u)$ is a non-separable complete metric space and
a
topological group.\\
$(3)$ $Aut_0(X,\B,\mu)$ is closed in $(Aut(X,\B,\mu), d_u)$;\\
$(4)$ $Aut(X,\B,\mu)$ (and $Aut_0(X,\B,\mu)$) has no outer group
automorphisms and has no closed normal subgroups with respect to the
uniform topology.\\
$(5)$ The topological space $(Aut(X,\B,\mu), d_u)$ is path-connected and
simple connected.\\
$(6)$ The topological space $Aut_0(X,\B,\mu)$ is contractible in the
uniform topology.
\end{theorem}

Some of these results are rather simple exercises which are known for a
long time. Statements (1) and (2) were firstly proved in \cite{Hal 1} for
the group $Aut_0(X,\B,\mu)$ and then they were used by many authors in
more general settings (see e.g. \cite{IT}, \cite{Ch-Kak}, \cite{B-Gol},
\cite{Alp-Pr 2}). The algebraic and topological structures of
$Aut(X,\B,\mu)$ (and $Aut_0(X,\B,\mu)$) were studied in \cite{Cha-Fr},
\cite{B-Gol}, \cite{E 1}, \cite{E 2}, and some other papers. The fact that
the space $\Aut$ is path-connected and even simple connected was proved in
\cite{B-Gol}. Notice also that the contractibility of $Aut_0(X,\B,\mu)$
was established in \cite{Kea}.

Let $U_T$ be the unitary operator $U_T$ acting on $L^1(X,\B,\mu)$ as
follows:
$$(
U_Tf)(x) = f(T^{-1}x)\sqrt{\frac{d\mu\circ T^{-1}}{d\mu}(x)}.
$$
$U_T$ is called  the operator associated to $T\in \Aut$.

The group $Aut(X,\B,\mu)$ endowed with weak topology has the following
topological properties.

\begin{theorem}[weak topology]\label{Aut-w}
$(1)$ $\Aut$ is a Polish space  and a
topological group with respect to the weak topology.\\
$(2)$ A sequence of automorphisms $(T_n)$ converges to $T$ in the weak
topology if and only if $(U_{T_n})$ converges to $U_T$ in the
strong (weak) operator topology.\\
$(3)$ The groups $Aut_0(X,\B,\mu)$ and $Aut(X,\B,\mu)$ are  contractible
in the weak topology.\\
$(4)$ $Aut_0(X,\B,\mu)$ is closed in $\Aut$ with respect to the weak
topology.\\
$(5)$ The set $\bigcup_{\nu\sim\mu} Aut_0(X, \B,\nu)$ is a meager
subset in $(\Aut, d_w)$.\\
\end{theorem}

Notice that statement (3) was proved in \cite{Kea} for the group
$Aut_0(X,\B,\mu)$ and then later it was generalized to $Aut(X,\B,\mu)$ in
\cite{Dan}. Statement (5), which asserts, in other words, that the set of
automorphisms of type $II$ is meager in $\Aut$, was established in
\cite{IT} (the first example of a non-singular automorphism which does not
in $\bigcup_{\nu\sim\mu} Aut_0(X, \B,\nu)$ was found by Ornstein). Observe
that since $(\Aut, d_w)$ is separable,  one can find various classes of
automorphisms which are dense in $\Aut$ in $d_w$. For example, assuming $X
= [0,1]$, the set of cyclic binary permutations on $[0,1]$ is dense in
$d_w$. In \cite{Pr 1}, the following result was proved: {\it The set of
pairs of elements from $Aut_0(X,\B,\mu) \times Aut_0(X,\B,\mu)$ which
generate a dense subgroup of $Aut_0(X,\B,\mu)$ form a dense $G_\delta$
subset.}

%%%%%%%%%%%%%%%%%%%%%%Rokhlin lemma%%%%%%%%%%%%%%%%%%%%%

\subsection{The Rokhlin lemma and Rokhlin property.} In \cite{Ro 1},
the following fundamental result (called lately the Rokhlin lemma) was
proved: {\it Let $T$ be an aperiodic finite measure preserving
automorphism of $(X,\B,\mu)$. Then given $n\in \N$ and $\e >0$, there
exists a set $A\in \B$ such that $A\cap T^iA = \emptyset,\ i=1,...,n-1$
and }
$$
\mu(X \setminus \bigcup_{i=0}^\infty T^iA) <\e.
$$

Based on this result, Rokhlin \cite{Ro 1} showed  that the set
$\per$ of periodic automorphisms is {\it dense} in
$Aut_0(X,\B,\mu)$ with respect to the uniform topology: {\it For
$T\in \ap, n\in \N$, and $\e >0$, there exists $P\in \per$ such
that $d_u(T,P) < \e+ 1/n$.}

The Rokhlin lemma lies in the basis of many areas of ergodic
theory. The key concept related to this result is a tower
constructed by an automorphism $T$: a collection of disjoint sets
$(A,TA,...,T^{n-1}A)$ is called a tower (or $T$-tower) of height
$n$ with base $A$.

The formulated above results can be generalized in several directions.
Firstly, one can prove that every non-singular automorphism is also
approximated by  periodic automorphisms in the uniform topology
\cite{Cha-Fr} (the Rokhlin lemma for non-singular automorphisms), that is
the set of periodic automorphisms, $\per$, is dense in $(\Aut, d_u)$.
Secondly, S.~Alpern \cite{Alp}, \cite{Alp-Pr 3} proved a multiple Rokhlin
tower theorem for measure preserving and non-singular aperiodic
automorphisms which shows that the space $X$ can be partitioned into a
finite or infinitely countable number of towers of prescribed measures.

\begin{theorem}[multiple Rokhlin tower lemma]\label{multitower}
For any $k\geq 2$, let $n_1,...,n_k$ be relatively prime positive integers
and let $p_1,...,p_k$ be positive  numbers such that $n_1p_1 +...+n_kp_k =
1$. Then for any aperiodic  $T\in Aut_0(X,\B,\mu)$ there exists sets $B_i
\in \B, i=1,...,k$,  with $\mu(B_i) = p_i$ and such that $\{T^jB_i :
i=1,...,k; j= 0,...,n_i-1\}$  is a partition of $X$.
\end{theorem}

In the next sections, we will discuss the Rokhlin lemma type results in
Borel and Cantor dynamics. It will be shown below that, like in ergodic
theory, the underlying space is partitioned into a finite or countable
collection of towers for these dynamics.

Let $G$ be a topological group. Following \cite{Gl-Ki}, we say that {\it
$G$ has the {\it Rokhlin property} if there exists some $g\in G$ such that
the conjugates of $g$ constitute a dense subset in $G$}.

\begin{theorem}[conjugacy lemma]\label{rokhl-prop} Let $T\in \Aut$ be an
aperiodic automorphism. Then:\\
$(1)$  $\{STS^{-1} : S \in \Aut\}$ is dense in $\ap$ with
 respect to the uniform topology;\\
$(2)$ $\{STS^{-1} : S \in \Aut\}$ is dense in $\Aut$  with respect to the
weak topology.
\end{theorem}

These results are due to Halmos \cite{Hal 2} for measure preserving
automorphisms. Theorem \ref{rokhl-prop} was generalized to the group $\Aut$
in \cite{Ch-Kak} and \cite{Fr}. Also the case of an infinite
$\sigma$-finite invariant measure was considered in \cite{Ch-Kak}. Note
that these statements are very useful for proving the density and
$G_\delta$-ness of various classes of automorphisms.

\subsection{\bf Classes of automorphisms.}
In this part, we will consider the following classes of automorphisms of a
measure space: periodic, aperiodic, ergodic, mixing and weakly mixing,
zero entropy, and of a finite rank.

As mentioned above, periodic automorphisms form a dense subset in
$\Aut$ with respect to the uniform and weak topologies:
\begin{equation}\label{dense-per}
  \overline{\per}^{d_u} = \overline{\per}^{d_w}= \Aut.
\end{equation}

Furthermore, the following two simple facts hold for the set of aperiodic
automorphisms: (1) {\it the set $\ap$ is a closed nowhere dense subset in
$\Aut$ with respect to the uniform topology; $(2)$ the set $\ap$ is dense
in $\Aut$ with respect to $d_w$.} Clearly, the first result is an immediate
consequence of the Rokhlin lemma. On the other hand, the second statement
follows from the Rokhlin property.

Recall that a set of a topological space is called  {\it generic} (or {\it
typical}) if it contains a dense $G_\delta$ subset.

In \cite{Hal 1}, Halmos answered the very important question about
genericity of ergodic transformations.

\begin{theorem}\label{erg-meas} The set of ergodic automorphisms is a
dense $G_\delta$ subset of $Aut_0(X,\B,\mu)$ with respect to the
weak topology and a nowhere dense subset in the uniform topology.
\end{theorem}

In the paper \cite{Ch-Kak}, this result was generalized to the group $\Aut$
and to the group $Aut_0(X,\B,\nu)$ where  $\nu$ is an infinite
$\sigma$-finite measure equivalent to $\mu$. Furthermore, the
Choksi-Kakutani result on genericity of ergodics was improved in
\cite{Ch-Haw-Pr} by showing that the ergodic type $III_1$ automorphisms
form a dense $G_\delta$ subset in $(\Aut, d_w)$.

In \cite{Ro 1}, it was proved  that finite measure preserving ergodic
automorphisms constitute a dense $G_\delta$ subset in $\ap$ with respect to
the uniform topology. It can be easily seen that ergodics are dense in
non-singular aperiodic automorphisms in $d_u$.

Let $K$ be a compact metric space. A Borel probability measure $\mu$ on $K$
is called an Oxtoby-Ulam measure (O-U measure) if it is non-atomic and
positive on every non-empty open set. If $K$ is a connected compact
manifold, then additionally the O-U measure $\mu$ must be zero on the
boundary $\partial K$. The following result, proved in \cite{O-U}, states
that a typical homeomorphism of $K$ is ergodic:

\begin{theorem}\label{erg} Let $\mu$ be an O-U measure on a compact
manifold $K$, $\dim K \geq 2$. Then the set of $\mu$-ergodic homeomorphisms
of $K$ is a dense $G_\delta$ subset in the group $Homeo_\mu(K)$ of all
$\mu$-preserving homeomorphisms of $K$ with respect to the topology of
uniform convergence. In particular, every compact manifold supports an
ergodic homeomorphism.
\end{theorem}

It is worth to mention one of the main results from \cite{Alp-Pr 2}.
Namely, it was proved that the natural embedding of $Homeo_\mu(K)$ into the
space $(Aut_0(K,\mu), d_w)$ preserves $G_\delta$-ness.

\begin{theorem} Let $\mu $ be an O-U measure on a compact connected
manifold $K$. Let $\mathcal V$ be a conjugate invariant, weak topology
$G_\delta$ subset subset of $Aut(K,\mu)$, which contains an aperiodic
automorphism. Then $\mathcal V\cap Homeo_\mu(K)$ is a dense $G_\delta$
subset of $Homeo_\mu(K)$ with respect to the topology of uniform
convergence.
\end{theorem}

The statistical (asymptotical)  properties of an ergodic automorphism are,
in general, rather poor. The dynamical systems generated by mixing and
weakly mixing automorphisms have more developed statistical properties. The
study of such automorphisms was begun in the papers \cite{Hal 2}, \cite{Hal
3}, \cite{Ro 1}, \cite{Ro 2}. Now there exist a large number of articles
devoted to mixing and weakly mixing automorphisms. We mention here only the
classical results.

\begin{theorem}\label{mix} $(1)$ The set of weakly mixing automorphisms
forms a dense $G_\delta$ subset in $Aut_0(X,\B,\mu)$ with respect to the
weak topology;\\
$(2)$ the set of mixing automorphisms is a meager subset of
$Aut_0(X,\B,\mu)$ with respect to the uniform topology.
\end{theorem}

Another crucial notion of ergodic theory is {\it entropy}
discovered by Kolmogorov and studied later in works by Rokhlin,
Sinai, Ornstein, and many others. In \cite{Ro 3}, the following
generic set was found:

\begin{theorem}\label{0-entr} Automorphisms of zero entropy form a
dense $G_\delta$ subset of $Aut_0(X,\B,\mu)$ with respect to the weak and
uniform topologies.
\end{theorem}

In the papers \cite{Kat-St 1} and \cite{Kat-St 2}, a new powerful method,
based on the Rokhlin lemma, was discovered. The authors called this method
a cyclic approximation by periodic automorphisms. Katok and Stepin showed
that the speed of approximation of an automorphism $T\in Aut_0(X,\B,\mu)$
by periodic ones determines ergodic, spectral, and mixing properties of
$T$. In particular, the following result was proved:

\begin{theorem} The set of automorphisms admitting cyclic approximation
with a fixed speed contains a dense $G_\delta$ subset in $Aut_0(X,\B,\mu)$
with respect to $d_w$.
\end{theorem}

It follows from this theorem that the set of automorphisms of rank 1 is
generic in $(Aut_0(X,\B,\mu), d_w)$.

In the space $Aut(X,\B, \nu)$, $\nu(X) = \infty$, one can also consider the
class of (infinite) measure preserving automorphisms. In this settings,
there exists a number of results on typical behavior of infinite measure
preserving automorphisms. We mention the following facts proved in
\cite{Kre} and \cite{Sa}, respectively.

\begin{theorem} $(1)$ Incompressible (conservative) automorphisms acting
on an infinite measure space form a dense $G_\delta$ subset in the
uniform and weak topologies.\\
$(2)$ Automorphisms with finite ergodic index form a meager set, but those
with infinite index are generic in  $(Aut_0(X,\B, \nu), d_w)$, $\nu(X) =
\infty$.
\end{theorem}

Recently, Choksi and Nadkarni \cite{Ch-N} proved that the set of
non-singular automorphisms with infinite ergodic index is a dense
$G_\delta$ subset in $\Aut$ with respect to the weak topology.

%%%%%%%%%%%%%%%%%%%%%%%%%%%Subgroups%%%%%%%%%%%%%%%%

\subsection{Subgroups of $\Aut$ defined by an automorphism.}
To every  automorphism $T$, we can associate the full group $[T]$ generated
by $T$, the normalizer $N[T]$, and the centralizer $C(T)$. Here we consider
topological properties of these groups. Notice that some of the results
formulated below hold for more complicated groups of automorphisms than
$\mathbb Z$.

\begin{theorem}\label{full-gr} $(1)$ The full group $[T]$ of any
automorphism $T\in \Aut$ is a closed nowhere dense subset of $\Aut$ in the
uniform topology;\\
$(2)$ the full group $[T]$ is a Polish group in the uniform topology;\\
$(3)$ the full group $[T]$ of any ergodic finite measure preserving
automorphism $T$ is a meager dense subset in $\Aut$ with respect
to the weak topology;\\
$(4)$ the full group $[T]$ of an ergodic automorphism $T$ is topologically
simple in $d_u$;\\
$(5)$ the full group $[T]$ is contractible in the uniform topology.
\end{theorem}

The first two results are  simple exercises. In statement (3), the fact
that $[T]$ is dense is also trivial, and the other part of the statement
was observed by A.~Kechris. The fact that for ergodic $T$ the full group
$[T]$ has no uniformly closed normal subgroups follows from \cite{Dye} (for
type $II$ automorphisms) and \cite{B-Gol} (for type $III$ automorphisms).
The last statement was proved in \cite{Dan}.

We note also that the notion of full group is extremely useful in the study
of problems related to the orbit equivalence theory of ergodic non-singular
automorphisms. As shown in \cite{G-P-S 2}, full groups (and topological
full groups) play a similar role in Cantor dynamics.

Recall that on the normalizer $N[T]$ can be endowed with a
topology generated by the metric $\rho$ \cite{Ha-Os}:
\begin{equation}\label{normal-topol}
\rho(R_1, R_2) = d_w(R_1,R_2) + \sum_{n\in \Z} 2^{-|n|}
\frac{d_u(R_1T^nR_1^{-1}, R_2T^nR_2^{-1})}{1 + d_u(R_1T^nR_1^{-1},
R_2T^nR_2^{-1} )}.
\end{equation}

The results formulated in the following theorem are taken from
\cite{Ha-Os} (the first two statements) and from \cite{Dan} (the
last two statements).

\begin{theorem}\label{normal} Let $T$ be an ergodic non-singular
automorphism. Then:\\
$(1)$ $N[T]$ is a Polish space with respect to $\rho$;\\
$(2)$ let $S\in Aut(X,\B,\nu), \nu(X) = \infty,$ and $\nu$ is
$S$-invariant; then $R\in \N[S]$ belongs to the closure of $[S]$ in $(N[S],
\rho)$ if and only if $\nu\circ R =\nu$;\\
$(3)$ the normalizer $N[T]$, where $T$ is of type $II$, is
contractible with respect to $\rho$;\\
$(4)$ if $T$ is an ergodic automorphism of type $III_\lambda$, then the
fundamental group $\pi_1(N[T]) = \Z$.
\end{theorem}

We finish this section with the weak closure theorem proved in
\cite{Ki} for rank 1 transformations.

\begin{theorem}\label{rank1} Let $T$ be a rank 1 automorphism from
$Aut_0(X,\B,\mu)$ and let $C(T)$ denote the centralizer of $T$.
Then $\overline{\{T^n :n\in \Z\}}^{d_w} \supset C(T)$.
\end{theorem}

Notice also that the centralizer $C(T)$ is closed in the weak topology.

%%%%%%%%%%%%%%%%%%%%%Section: Borel dynamics%%%%%%%%%%%%%%%

\section{Topologies in Borel dynamics}
In this section, we define and study  topologies on $\aut$ which are
analogous to the weak and uniform topologies on $\Aut$. There exist many
papers devoted to Borel dynamics but, as far as we know, topologies on
$\aut$ were first studied in our papers  \cite{B-K 1}, \cite{B-D-K 1}, and
 \cite{B-M}.

\subsection{The uniform and weak topologies on $\aut$.}
We first define two topologies on $\aut$ and then consider the topological
properties of the space $\aut$.

Given $S,T \in \Aut$, denote by $E(S,T)=\{ x\in X \ |\ Tx\ne Sx\}
\cup \{x\in X \ |\ T^{-1}x\ne S^{-1}x\}$.

\begin{definition}\footnote{In fact,
we defined in \cite{B-D-K 1} more topologies on $\aut$ but our main results
concern mostly $\tau$ and $p$, so that we do not consider the remaining
topologies in this paper.} \label{DefTopologies}  The  topologies $\tau$
and $p$ on $\aut$ are defined by the bases of neighborhoods $\mathcal U$
and $\mathcal W$, respectively. They are: ${\mathcal U} = \{U(T;
\mu_1,...,\mu_n; \varepsilon)\}$,  $\mathcal W = \{W(T; F_1,...,F_k)\}$,
where
\begin{equation}\label{DefUnifTop}
U(T; \mu_1,...,\mu_n; \varepsilon)= \{ S\in \aut\ |\
\mu_i(E(S,T))<\varepsilon ,\ i=1,...,n\},
\end{equation}
\begin{equation}\label{DefWeakTop}
W(T; F_1,...,F_k) = \{S\in \aut\ |\ SF_i = TF_i, \ i= 1,...,k\},
\end{equation}
In all the above definitions $ T\in \aut $, $ \mu_1,...,\mu_n \in \meas,\
F_1,...,F_k \in \B$,  and $\varepsilon >0 $.
\end{definition}

We call $\tau$ and $p$ the {\it uniform and weak topologies}, respectively.
Our motivation is as follows. For a standard Borel space, we do not have a
fixed Borel measure on the underlying space. Therefore, if we want to
extend  definitions (\ref{d_u}) and (\ref{d_w}) to $\aut$, we have to take
into account the set ${\mathcal M}_1(X)$ of {\it all} Borel probability
measures on $(X,\B)$.  To generalize the definition of $d_w$, we observe
that if the symmetric difference of two Borel sets is arbitrarily small
with respect to {\it any} $\mu\in {\mathcal M}_1(X)$, then these sets must
coincide.

Notice that if the set $E_0(S,T) = \{x\in X : Sx\neq Tx\}$ were used in
(\ref{DefUnifTop}) instead of $E(S,T)$, then we would obtain the topology
equivalent to $\tau$.

Observe that the uniform topology $\tau$ is equivalent to $\tau'$ which can
be defined on $\aut$ similarly to the metric $d'_u$ (\ref{d'_u}).

It is natural to consider two more topologies $\tau_0$ and $p_0$ on $\aut$
(they are some modifications of $\tau$ and $p$) by considering only
continuous measures and uncountable Borel sets in Definition
\ref{DefTopologies}.

\begin{definition}\label{ContTopolBorel} The topologies $\tau_0$ and
$p_0$ on $\Aut$ are defined by the bases of neighborhoods ${\mathcal U}_0
= \{U_0(T;\nu_1,...,\nu_n;\e)\}$ and ${\mathcal W}_0 =
\{W_0(T;A_1,...,A_n)\}$, respectively, where $U_0(T;\nu_1,...,\nu_n;\e)$
and $W_0(T;A_1,...,A_n)$ are defined as in (\ref{DefUnifTop}) and
(\ref{DefWeakTop}) with continuous measures $\nu_i \in \meas$ and
uncountable Borel sets $A_i$ for $i=1,...,n$.
\end{definition}

Clearly, $\tau$ and $p$ are not weaker than $\tau_0$ and $p_0$,
respectively.

Given an automorphism $T$ of $\bs$, we can associate a linear
unitary operator $L_T$ on the Banach space $B(X)$ of all bounded
Borel functions by $(L_Tf)(x) = f(T^{-1}x)$. Let $\tilde p$ be the
topology induced on $\aut$ by the strong operator topology on
bounded linear operators of $B(X)$. It was shown in \cite{B-D-K 1}
that the topologies $p$ and $\tilde p$ are equivalent. This fact
corresponds to the well known result in ergodic theory (see
Theorem \ref{Aut-w}) and is another justification of the name {\it
\textquotedblleft weak topology"} which is used to refer to $p$.

The following theorem proved in \cite{B-D-K 1} reveals some
topological properties of $\aut$ with respect to the topologies we
have defined.

\begin{theorem}\label{UnifBor} $(1)$ $\aut$ is a Hausdorff topological
group with respect to the topologies $\tau,\tau_0,p$, and $p_0$;\\
$(2)$ the sets $W(T;F_1,\ldots;F_n)$ are closed in $\aut$ with respect to
$\tau$ and clopen with respect to $p$, that is $(\aut, p)$ is a
0-dimensional topological group;\\
$(3)$ $\aut$ is a complete non-separable group with respect to $\tau$ and
$p$ in the sense that every Cauchy sequence of Borel
automorphisms converges to a Borel automorphism;\\
$(4)$ the group $\aut$ is totally disconnected in the uniform topology
$\tau$;\\
$(5)$ the group $\aut$ is path-connected in the topology $\tau_0$;\\
$(6)$ the topology $p$ is equivalent to $p_0$;\\
$(7)$ the topology $\tau$ is not comparable with $p$.
\end{theorem}

Let us give a few comments on the results from Theorem \ref{UnifBor}.
Firstly, observe that, in contrast to the ergodic theory, the uniform and
weak topologies are not comparable on $\aut$. Secondly, to see that that
$(\aut, \tau)$ is totally disconnected, it suffices to notice the following
facts: (i) the set of all automorphisms which have a fixed period at a
fixed point is clopen in $\tau$ (see Theorem \ref{periodic}(3) below) and
(ii) $Ctbl(X)$ is the unique normal subgroup of $\aut$ which is not
connected in its turn. On the other hand, the fact that the topological
space $(\aut, \tau_0)$ is path-connected (proved in \cite{B-M}) well
corresponds to the results from \cite{B-Gol}, \cite{Dan}, \cite{Kea}
mentioned in Theorem \ref{Aut-un}.

It is obvious that the topologies on $\aut$ are not defined by convergent
sequences, however, it is useful for many applications to know criteria of
convergence.

\begin{remark}{\rm (1) $(T_n)$ converges to $S$ in $\tau$ if and
only if $\forall x \in X\ \exists n(x)\in \N$ such that $\forall n
>n(x), \ T_nx = Sx$.

(2) $(T_n)$ converges to $S$ in $p$ if and only if for any Borel set $F$,
$T_nF = TF$ for all sufficiently large $n$. In particular, $F$ can be a
point from $X$. Therefore, we see that $p$-convergence implies
$\tau$-convergence.

(3) Observe that the criterion of convergence (1) does not hold for the
topology $\tau_0$.}
\end{remark}

\subsection{The quotient group $\widehat{Aut}(X,\B)$ and quotient
topologies.}

Identifying Borel automorphisms  from $\aut$ which are different on an at
most countable set, we obtain elements from the quotient group
$\widehat{Aut}(X,\B) = \aut/Ctbl(X)$. Observe that $Ctbl(X)$ is a closed
normal subgroup with respect to the topologies $\tau,\ p$, and $\tau_0$
(see \cite{B-M} and \cite{B-D-K 1}). Denote by $\widehat\tau, \widehat
{\tau_0}$ and $\widehat p$ the quotient topologies induced on
$\widehat{Aut}(X,\B)$ by $\tau, \tau_0$ and $p$, respectively. Notice that
the bases of the topologies $\widehat\tau$ and $\widehat p$ form by the
sets $\widehat U(T;\mu_1,...,\mu_n;\e) = U(T; \mu_1,...,\mu_n;\e)Ctbl(X)$
and $\widehat W(T;F_1,...,F_m) = W(T;F_1,...,F_m)Ctbl(X)$, respectively,
with continuous measures $\mu_i$ and uncountable Borel sets $F_j$. The
class of automorphisms equivalent to a Borel automorphism $T$ we again
denote by the same symbol $T$ and write $T\in \widehat{Aut}(X,\B)$. This
identification corresponds to the situation in ergodic theory when two
automorphisms are identified if they are different on a set of measure 0
(mod-0-convention).

The following theorem clarifies some topological and algebraic properties
of $\widehat{Aut}\bs$.

\begin{theorem}
$(1)$ $\widehat{Aut}\bs$ is a Hausdorff topological group with respect to
the quotient topologies $\widehat\tau$ and $\widehat p$; the topologies
$\widehat {\tau_0}$ and  $\widehat\tau$ are  equivalent on
$\widehat{Aut}\bs$;\\
$(2)$ $\widehat{Aut}\bs $ is an algebraically simple group;\\
$(3)$ the group $(\widehat{Aut}(X,\B),\widehat\tau)$  is path-connected.
\end{theorem}

The second statement was proved in \cite{Sh} and the other results were
obtained in \cite{B-M}. Observe that the identification of automorphisms
different on an at most countable set improves considerably the topological
structure of the quotient group. Apparently, one can conjecture the
contractibility of $\widehat{Aut}\bs$ in $\widehat\tau$. This would be
parallel to the known result in the ergodic theory. We believe that the
space $(\widehat{Aut}\bs, \widehat\tau)$ has the topological and algebraic
structures similar to those of $(\Aut, d_u)$.

\subsection{Approximation results in $\aut$.} In this section, we
consider several natural classes of automorphisms and find their closures
in the topologies we have defined above.
\medskip

We first consider the sets of periodic and aperiodic Borel automorphisms.
Let $T\in \Aut$. Then the space $X$ can be partitioned into a disjoint
union of Borel $T$-invariant sets $X_1, X_2,...,X_\infty$ where $X_n$ is
the set of points with period $n$, and $X_\infty$ is the set where $T$ is
aperiodic. Denote by $\per_n(x)$ the set of all automorphisms which have
period $n$ at $x$. By definition, $T$ is  in $\per_n$, the set of all Borel
automorphisms of period $n$, if $X_n = X$. In other words,
\begin{equation}\label{per_n}
\per_n =\bigcap_{x\in X} \per_n(x).
\end{equation}
We say that $T\in \per_0$ if there exists $N\in \N$ such that
$P^Nx=x,\ x\in X$.\smallskip

The following two theorems describe some topological properties of the sets
of periodic and aperiodic automorphisms in $\aut$ (cf. (\ref{dense-per})).

\begin{theorem}[periodic automorphisms]\label{periodic}
$(1)$ The set  $\per_0$ is dense in $(\aut,\tau)$.\\
$(2)$ $\per$ is a closed nowhere dense subset in $\aut$ with respect to
$p$.\\
$(3)$ For any $n\in \N$, the set $\per_n(x)\ (x\in X)$ is clopen with
respect to the topologies $\tau$ and $p$.\\
$(4)$  $\overline{\per_0}^\tau = \overline{\per}^\tau$.\\
$(5)$ $\per_n\ (n\in \N)$ is closed with respect to $\tau$ and $p$.\\
\end{theorem}

Notice that statement (2) was proved in \cite{B-M} and all other results
are taken from \cite{B-D-K 1}. It is curious enough that the set of
periodic automorphisms is rare in the topology $p$.

\begin{theorem}[aperiodic automorphisms]\label{TopPropAper}  The set $\ap$
is closed and nowhere dense in $\aut$ with respect to $\tau$ and $p$.
\end{theorem}

{\it Proof.}  We prove here that $\ap$ has no interior points in the
topology $p$ (all other statements were obtained in \cite{B-D-K 1}). Let
$T\in\ap $ and let $W=W(T;F_1,\ldots,F_n)$ be a $p$-neighborhood of $T$,
where $(F_1,\ldots, F_n)$ is a partition of $X$. We need to show that $W$
contains a non-aperiodic automorphism. To do this, it suffices to change
$T$ at some points and construct an automorphism $P\in W$ which has a
periodic orbit. Fix a point $x\in X$ and consider the $T$-orbit of $x$.
Then we can find two numbers $i,j\in\mathbb Z$ and two sets (say $F_{n_1}$
and $F_{n_2}$) such that $j\geq i+2$ and $T^ix, T^jx\in F_{n_1}$ and
$T^{i+1}x, T^{j+1}x\in F_{n_2}$ (possibly, $n_1 = n_2$). Define the
automorphism $P$ as follows: $PT^ix=T^{j+1}x$, $PT^jx=T^{i+1}x$, and
$Px=Tx$ elsewhere. It is not hard to check that $P\in W$ and, clearly,
$P\notin \ap$. \hfill$\square$
\medskip

The next lemma is one of the main tools in the  study of Borel
automorphisms. The proof can be found in \cite{Be-Ke} or \cite{N}.
We formulate the lemma as in \cite{B-D-K 1}.

\begin{lemma}\label{markers-lemma} Let $T\in \Aut$ be an aperiodic Borel
automorphism of a standard Borel space $\bs$. Then there exists a sequence
$(A_n)$ of Borel sets such that\\
(i) $X=A_0 \supset A_1 \supset A_2\supset \cdots,$\\
(ii) $\bigcap_n A_n =\emptyset,$\\
(iii) $A_n$ and $X\setminus A_n$ are
complete $T$-sections, $n\in \N$,\\
(iv) for $n\in \N$, every point in
$A_n$ is recurrent,\\
(v) for $n\in \N$,\ $A_n\cap T^i(A_n) =\emptyset,\ i=1,...,n-1$,\\
(vi) for $n\in \N$, the base  of every non-empty $T$-tower built by the
function of the first return to $A_n$ is an uncountable Borel set.
\end{lemma}

A sequence of Borel sets $(A_n)$ satisfying  Lemma \ref{markers-lemma} is
called a {\it vanishing sequence of markers}.
\smallskip

Recall that a Borel partition $\Xi=\{\xi_i : i\in I\}$\ ($|I| \leq
\aleph_0$) of $X$ is called a {\it K-R (Kakutani-Rokhlin) partition} for an
automorphism $T\in\aut$ if all the $\xi_i$'s are $T$-towers,
$\xi_i=\{B_i,TB_i,\ldots,T^{m-1}B_i\}$, for some $B_i\in\B$. Set
$B(\Xi)=\bigcup_{i\in I}B_i$.

The crucial importance of Lemma \ref{markers-lemma} may be partially
explained by the following fact. It turns out that, having a vanishing
sequence of markers, one can construct a sequence of Kakutani-Rokhlin
partitions satisfying the following properties.

\begin{proposition}[K-R partitions] \label{K-R partitions} Let $T$
be an aperiodic Borel automorphism of $\bs$. There exists a sequence of
K-R partitions $\{\Xi_n\}$ of $X$ whose elements generate the
$\sigma$-algebra $\B$ and such that: (i) $\bigcap_n B(\Xi_n)=\emptyset$;
(ii) $\Xi_n$ refines $\Xi_{n+1}$; and (iii) if $h(\xi_i)$ is the height of
$\xi_i$, then $h_n=\min\{h(\xi) : \xi\in \Xi_n\}\to\infty$ as
$n\to\infty$.
\end{proposition}

Based on this result, we introduced in \cite{B-D-K 1} the concept of
Bratteli diagrams in the settings of Borel dynamics. It was shown that one
can associate an ordered Bratteli diagram to every Borel aperiodic
automorphism $T$ such that $T$ becomes isomorphic to the Vershik map acting
on the space of infinite paths of the Bratteli diagram. Notice that the set
of vertices at each level of the diagram (i.e., the number of towers in the
K-R partitions)  can be chosen finite, and the ordered diagram
corresponding to $T$ has no maximal and minimal paths. It is worthy to
recall that the notion of Bratteli diagrams was originally introduced in
the theory of approximately finite $C^*$-algebras and then extensively
studied by many authors in the setting of Cantor minimal systems (see, e.g.
\cite{H-P-S}, \cite{G-P-S 1}, \cite{Ma}).
\smallskip

We conclude this section with the following theorem \cite{B-D-K
1}.

\begin{theorem} Let $T\in \aut$. Then:
$(1)$ the full group $[T] \ (T\in\aut)$ is closed and
nowhere dense in $\aut$ with respect to the topologies $\tau$ and $p$;\\
$(2)$ $\per_0\cap [T]$ is $\tau$-dense in $[T]$ for each aperiodic $T$.
\end{theorem}

%%%%%%%%%%%%%%%%%%%%Rokhlin%%%%%%%%%%%%%%%%%%%%%

\subsection{\bf The Rokhlin lemma and Rokhlin property.}
The first result on approximation of any aperiodic Borel automorphism by
periodic ones  in the context of Borel dynamics was proved by Weiss in
\cite{W}. Following \cite{W}, we call a Borel set $C$ {\it completely
positive} with respect to an aperiodic automorphism $T$ if $\mu(C) > 0$ for
any $T$-nonsingular measure $\mu$. A Borel subset $A$ is called {\it
wandering} for $T$ if $T^iA\cap A = \emptyset, i\neq 0$.

\begin{theorem} If $T$ is aperiodic Borel automorphism and $C$  is a Borel
set such that  $X \setminus C$ is completely positive, then for any prime
$p$ there exists a Borel set $B$ with (i) $B,TB,...,T^{p-1}B$ pairwise
disjoint, and (ii) the set $C \setminus (\bigcup_{j=0}^{p-1} T^jB)$ is
wandering.
\end{theorem}

To see how this theorem is related to the classical Rokhlin lemma, we
recall the following fact. Let $\mathcal W(T)$ denote the $\sigma$-ideal of
wandering sets with respect to $T$, and let $\mathcal Q(T)$ be the set of
$T$-nonsingular measures on $\B$. Then
$$
\mathcal W(T) = \bigcap_{\mu\in \mathcal Q(T)} \mathcal N(T),
$$
where $\mathcal N(T)$ is the $\sigma$-ideal of $\mu$-null sets \cite{W}.

Another approach to the problem of periodic approximation was applied by
Nadkarni in \cite{N}. In our notation, Nadkarni showed that {\it for given
$T\in \ap$, there exists a sequence of periodic automorphisms $(P_n)$ of $
(X,\B)$ such that $P_n \stackrel{\tau}{\longrightarrow} T,\ n\to \infty$.
Moreover, it the $P_n$'s can be taken from $[T]$}. It is clear now that
this result is a direct consequence of Proposition \ref{K-R partitions}.

It turns out that the Rokhlin lemma can be also proved for any
measure which is not related to a given aperiodic automorphism $T$
(see \cite{B-D-K 1} for details).

\begin{theorem}[Rokhlin lemma]\label{Rokhl-lemma}
Let $m\in \N$ and let $T$ be an aperiodic Borel automorphism of $\bs$.
Then for any $\e>0$ and any measures $\mu_1,...,\mu_p$ from $ \meas$ there
exists a Borel subset $F$ in $X$ such that $F, TF,..., T^{m-1}F$ are
pairwise disjoint and
$$
\mu_i(F\cup TF\cup\cdots \cup T^{m-1}F) > 1-\e,\ \ i=1,...,p.
$$
\end{theorem}

To emphasize the obvious similarity between measurable dynamics and Borel
dynamics, we give also two results on the Rokhlin property which were
proved in \cite{B-M}.

\begin{theorem}[Rokhlin property]\label{RokhlPropBor}
$(1)$ The action of $\widehat{Aut}\bs$ on itself by conjugation is
topologically transitive with respect to the topology $\widehat p$.
Moreover, $\{T^{-1}ST : T\in \widehat{Aut}\bs\}$ is dense in
$(\widehat{Aut}\bs,\widehat p)$ for any $S\in\mathcal Sm\cap\ap$.\\
$(2)$ Let $T$ be an aperiodic automorphism of $\bs$, then $\{S^{-1}TS :
S\in \aut\}$ is dense in $(\ap,\tau)$.
\end{theorem}

We should also refer to the recent article \cite{Ke-Ros} where the Rokhlin
property is considered in the context of automorphism groups of countable
structures.

%
%%%%%%%%%%%%%%%%%%%%%%%%%%%%%%   SMOOTH AUTOMORPHISMS
%

\subsection{\bf Smooth automorphisms.} We consider here the topological
properties of the class $\mathcal Sm$ of smooth automorphisms of a standard
Borel space. This class is a natural extension of the class of periodic
automorphisms.

We first note that the set $\mathcal Sm$ is dense in $(\aut,\tau)$. It
follows from the fact that $\per$ is dense in $(\aut,\tau)$ and
$\per\subset \mathcal Sm$. On the other hand, by Theorem \ref{TopPropAper}
the set $\mathcal Sm\cap \ap$ is not dense in $\aut$ with respect to
$\tau$ and $p$.

As shown in Theorem \ref{periodic}, the set of periodic automorphisms is
not dense in $\aut$ with respect to $p$. Kechris conjectured that to obtain
a dense subset in $(\aut, p)$, one needs to take the set of smooth
automorphisms. The following theorem was proved in \cite{B-M}.

\begin{theorem} $(1)$ $\overline{\mathcal Sm}^p=\aut$. \\
$(2)$  $\mathcal Sm\cap \ap$ is dense in $(\widehat{Aut}\bs,\widehat
p)$.\\
$(3)$ $\mathcal Sm\cap \ap$ is dense in $\ap$ with respect to $p$.
\end{theorem}

\subsection{\bf The topology of uniform convergence.} Let $(X,d)$ be a
compact metric space. In this case, we can consider the group $Aut(X,\B)$
of Borel automorphisms and its subgroup $Homeo(X)$ of homeomorphisms of
$X$. Define for $S,T\in \aut$ the topology of uniform convergence generated
by the metric
\begin{equation}\label{metricD}
D(S,T) = \sup_{x\in X} d(Sx,Tx) + \sup_{x\in X} d(S^{-1}x,T^{-1}x).
\end{equation}
Then $(\aut, D)$ is a complete metric space and $Homeo(X)$ is $D$-closed in
$\aut$. The case when $X (= \Om)$ is a Cantor set is considered in Section
4 (see \cite{B-D-K 2}, \cite{B-D-M} for details). It is not hard to see
that in Cantor dynamics the topology on $\h$ generated by $D$ is equivalent
to the topology $p$ defined by clopen sets only. But in Borel dynamics the
topologies $D$ and $p$ are different. Hence, it is interesting to find out
which topological properties  are preserved under embedding of $\h$ into
$Aut(\Om,\B)$ endowed with the topologies $D$ and $p$. We note also that
the topology generated by $D$ on $\aut$ depends, in general, on the
topological space $(X,d)$. Nevertheless, we think it is worthy to study the
topological properties of $\aut$ and $Homeo(X)$ for a fixed compact (or
Cantor) metric space  $X$ because one can compare in this case these
properties for the both groups.

The first problem we discuss for the metric $D$ is to find the closure of
aperiodic Borel automorphisms of $X$ with respect to $D$.

Let $T$ be a Borel automorphism of $X$. Then $X$ is decomposed into the
canonical $T$-invariant partition $(Y_1,Y_2,...,Y_\infty)$ where $T$ has
period $n$ on $Y_n$ and $T$ is aperiodic on $Y_\infty$ (see above {\bf
3.3}). We call $T$ {\it regular} if all the sets $Y_i,\ 1\leq i < \infty$,
are uncountable.

\begin{proposition}\label{regularBorel} {\rm \cite{B-D-K 1}}
Suppose that $T\in Aut(X,\B)$ is regular. Then for any $\e> 0$ there exists
$S\in \ap$ such that $D(S,T) < \e$.
\end{proposition}

Since $Ctbl(X)$ is closed in $D$, we immediately obtain

\begin{corollary} \label{} {\rm \cite{B-D-K 1}} The set of aperiodic
automorphisms from $\widehat{Aut}(X,\B)$  is dense with respect to the
quotient topology $\widehat D$ of $D$.
\end{corollary}

To answer the question when a non-regular automorphism $T\in \aut$ belongs
to $\overline{\ap}^D$, we need the following definition. We say that an
automorphism $T$ is {\it semicontinuous} at $x\in X$ if for any $\e
>0$ there exists $z\neq x$ such that $d(x,z) < \e$ and $d(Tx,Tz) < \e$.

\begin{theorem}\label{nonregularBorel} {\rm \cite{B-D-K 1}} Let $T$ be a
non-regular Borel automorphism from $\aut$ and let $(Y_1,Y_2,...,Y_\infty)$
be the canonical partition associated to $T$. Denote by $Y_0$ the set
$\bigcup_{i\in I}Y_i$ such that each $Y_i,\ i\in I,$ is an at most
countable set, $I \subset \N$. Then $T\in \overline{\ap}^D$ if and only if
for every $x\in Y_0$ there exists $y\in Orb_T(x)$ such that $T$ is
semicontinuous at $y$.
\end{theorem}

To finish this section we present two more results on closures of
some subsets of $\aut$ proved in \cite{B-D-K 1}.

\begin{theorem}\label{incompres} The set
$\mathcal Inc$ is a closed nowhere dense subset of $\ap$ with
respect to the topology $p$.
\end{theorem}

Say that $ T \in \ap \mod(Ctbl)$ if $T$ is aperiodic everywhere except an
at most countable set.

\begin{theorem}\label{odometers} $(1)$\ $\overline{{\mathcal O}d}^\tau =
\ap$ and $\overline{{\mathcal O}d}^{\tau_0} = \ap \mod(Ctbl)$.\\
$(2)$\ $\overline{{\mathcal O}d}^p \subset \mathcal Inc$.\\
$(3)$\  $\overline{{\mathcal S}m}^D \supset \overline{{\mathcal O}d}^D
\supset \ap$  assuming that $(X,d)$ is a compact metric space.
\end{theorem}

We note that full groups in the context of Borel dynamics were studied in
many papers. We mention here  \cite{J-Ke-L}, \cite{Fel-Moo}, \cite{Mer}.

%%%%%%%%%%%%%%%%%%%%%%%%%%%%%%%%%%%% CANTOR DYNAMICS %%%%%%%%%%%%%%%%%%%

\section{Topologies in Cantor dynamics.}

In this section, we consider the uniform and weak topologies in the
settings of Cantor dynamics.

%
% %%%%%%%%%%%%%%%%%%%%%%%   TOPOLOGICAL PROPERTIES OF $homeo(x$
%

\subsection{The weak and uniform topologies on $\h$.}
We define two topologies on $\h$ which are similar to the uniform and weak
topologies used in Borel dynamics. Our primary goal is to study the global
topological properties of $\h$ and compare these properties with  those
known for the groups $\aut$ and $\Aut$.

\begin{definition}{\rm
(i) The uniform topology $\tau$ on $\h$ is defined as the relative
topology on $\h$ induced from $(Aut(\Om,\B), \tau)$. The base of
neighborhoods is formed by
\begin{equation}\label{DefTau}
U(T; \mu_1,...,\mu_n; \e)= \{S\in \h)\ |\ \mu_i(E(S,T))<\e ,\ i=1,...,n\}.
\end{equation}
(ii) The topology $p$ is defined on $\h$ by the base of neighborhoods
\begin{equation}\label{DefWeak}
W(T; F_1,...,F_k) = \{S\in \h\ |\ SF_i = TF_i, \ i= 1,...,k\}.
\end{equation}
All measures in (\ref{DefTau}) are taken from $\meas$, and all the sets in
(\ref{DefWeak}) are assumed to be  clopen.}
\end{definition}

Apparently, the most natural topology on the group $\h$ is the topology $D$
of uniform convergence which is generated by the metric $D(T,S)=\sup_{x\in
\Om}d(Tx,Sx)+\sup_{x\in \Om}d(T^{-1}x,S^{-1}x)$ ($d$ is a metric on $\Om$
compatible with the topology). Recall that this topology can be induced
from $Aut(\Omega,\B)$ (\ref{metricD}).

It is not hard to show that on the group $\h$ the topologies $p$ and $D$
are equivalent and, as in the Borel case, the topologies $\tau$ and $p$ are
not comparable (see details in \cite{B-D-K 2}).

Now we formulate several statements concerning topological and algebraic
properties of $\h$.

\begin{theorem}\label{propertiesHomeo} $(1)$ The group $\h$ is simple.\\
$(2)$ $\h$ is a Hausdorff topological group with respect to the topologies
$\tau$ and $p$.\\
$(3)$ $(\h,p)$ is a $0$-dimensional perfect Polish space and
non-Archimedean topological group.\\
$(4)$ The topological space $(\h,\tau)$ is totally disconnected.\\
$(5)$ The group $\h$ is dense and non-closed in $(Aut(\Om,\B),\tau)$.
\end{theorem}

Let us make a few comments on these results. The first statement is rather
old and was proved in \cite{An}. The third result follows from the fact
that the neighborhoods in (\ref{DefWeak}) are closed, and the neighborhoods
of the identity are subgroups of $\h$. To see that $\h$ is totally
disconnected in $\tau$, we observe the following facts: (i) by (1), $\h$
has no normal subgroups, and (ii) there are nontrivial clopen subsets in
$\h$ (for example, the set $\{S \in \h : Sx =x\}$ is clopen with respect to
$\tau$). The last statement can be reformulated as follows: every Borel
automorphism of $\Om$ is approximated in $\tau$ by a homeomorphism. This is
a form of Luzin's theorem. A similar approach can be also found  in
\cite{Alp-Pr 1} and \cite{Go-Ni-Wa}.

Recall that $T_n\stackrel{\tau}{\longrightarrow}S \iff \forall x \in \Om\
\exists n(x)\in \N$ such that $\forall n >n(x), \ T_nx = Sx$. On the other
hand, $T_n\stackrel{p} {\longrightarrow} S \iff \forall {\rm clopen}\ F \
\exists n(F)$ such that $\forall n > n(F),\ T_nF = SF$, and $p$-convergence
is equivalent to the uniform convergence, i.e., in the metric $D$.

\subsection{The Rokhlin lemma and Rokhlin property.} We will discuss in
this section the problem of approximation of aperiodic homeomorphisms by
periodic ones. Recall that $\per$ denote the set of (pointwise) periodic
homeomorphisms and $\per_0$ is the set of periodic homeomorphisms with
finite period.

We begin with the following result proved in \cite{B-D-M} which
asserts the existence of K-R partitions for aperiodic
homeomorphisms.

\begin{theorem}\label{towers} Let $T$ be an aperiodic homeomorphism of
$\Om$. Given a positive integer $m \geq 2$, there exists a partition $\Xi
= (\xi_1,...,\xi_k)$ of $\Om$ (a K-R partition) into a finite number of
$T$-towers $\xi_i$ such that the height $h(\xi_i)$ of each tower is at
least $m$.
\end{theorem}

The idea of the proof is the following. Given $m\in\N$, find a clopen
partition of $\Om=\bigcup_{i=1}^qU_i$ such that $T^i U_j\cap U_j=\emptyset$
for $i=0,\ldots,m-1$. Set $A_1=U_1$ and
$A_i=U_i\setminus(\bigcup_{j=-m}^mT^j(A_1\cup\ldots\cup A_{i-1}))$ for
$i=2,\ldots,q$. Notice that the set $A = A_1\cup\ldots\cup A_q$ is clopen
and  meets every $T$-orbit at least once. Therefore for every $x\in A$
there is $l=l(x)>0$ such that $T^lx\in A$. The function of the first return
defines a K-R partition $\Xi$ of $\Om$ with base $A$ such that the height
of every $T$-tower from $\Xi$ is at least $m$.

Based on this result, we prove a topological version of the Rokhlin lemma
(the details are in \cite{B-D-M}) which contains the following three
results concerning approximation of aperiodic homeomorphisms by periodic
ones.

\begin{theorem}\label{RokhlinLemmaCantor} {\rm (Rokhlin lemma)}
Let $T$ be an aperiodic homeomorphism of $\Om$.\\
$(1)$ For $n\in \N$, $\e>0$, and any measures $\mu_1,...,\mu_k$ from
$\meas$, there exists a clopen subset $F$ in $\Om$ such that $F, TF,...,
T^{n-1}F$ are pairwise disjoint and
$$
\mu_i(F\cup TF\cup\cdots \cup T^{n-1}F) > 1-\e,\ \ i=1,...,k.
$$
$(2)$ There exists a sequence of periodic homeomorphisms $\{P_n\}\subset
Per_0\cap [[T]]$ such that $P_n\stackrel{\tau}{\longrightarrow} T$\
($[[T]]$ is the topological full group).\\
$(3)$ Given $n\in \N$ and $\e > 0$, there exists a K-R partition $\Xi =
(\xi_1,...,\xi_q)$ of $\Om$ as in Theorem \ref{towers} such that for
$i=1,\ldots,k$
$$
\mu_i(\bigcup_{s=1}^q\bigcup_{j=0}^{h(\xi_s)-n})>1-\e.
$$
\end{theorem}

To see that (3) holds, choose $m$ sufficiently large and repeat the
construction used in Theorem \ref{towers} to find some $i$ from
$(1,...,m-n-1)$ such that $\mu_j(T^iA\cup \ldots \cup T^{i+n}A)<\e$ for all
$j=1,\ldots,k$. To obtain the K-R partition $\Xi$, apply again the above
construction to the set $T^{i+n}A$ taken as a base of $\Xi$. The other
statements are consequences of this result.
\smallskip

\noindent{\bf Remark.} Observe that Theorem \ref{RokhlinLemmaCantor} allows
one to construct a Bratteli diagram associated to an aperiodic
homeomorphism. Moreover, one can show that every aperiodic homeomorphism of
a Cantor set can be represented as the Vershik map acting on the path space
of an ordered Bratteli diagram (see \cite{M} for more details).
\medskip

It turns out that, as in measurable and Borel dynamics, the Rokhlin
property also holds in the group $\h$.

\begin{theorem} $(1)$ Let $T$ be an aperiodic homeomorphism of $\Om$.
Then the set $\{S^{-1}TS : S\in \h\}$ is dense in $\ap$ with respect to
the uniform topology $\tau$.

$(2)$ There exists a dense $G_\delta$ subset $\mathcal E$ of $(\h, p)$ such
that $\{S^{-1}TS : S\in\h\}$ is $p$-dense in $\h$ for every $T\in\mathcal
E$.
\end{theorem}

The first assertion was obtained in \cite{B-D-M} and the second one was
proved in \cite{Gl-W 2}. Observe that it is a nontrivial problem to point
out explicitly a homeomorphism from the set $\mathcal E$. Indeed, if, say,
$T$ is minimal then the closure of the conjugacy class of $T$ does not
leave the set $\overline{\mathcal Min}^p$. But the latter is a nowhere
dense subset of $\h$ in the topology $p$ (see below). We do not know
whether there is an aperiodic homeomorphism such that its conjugates are
dense in $(\h, p)$.

The Rokhlin property was also studied for some connected compact spaces. In
\cite{Gl-W 2}, the Rokhlin property was also established for any finite
dimensional cube and the Hilbert cube. The local Rokhlin property was
proved for circle homeomorphisms in \cite{A-Hu-Ken}.

\subsection{Approximation results in the group $\h$.}

In this section, we consider various natural subsets of $\h$ and find
their closures. As mentioned above, $\h$ is not closed in $\aut$ in the
uniform topology $\tau$, therefore the $\tau$-closure of a subset
$Y\subset \h$ does not belong to $\h$, in general. For convenience, we
will  use the following convention: $\overline{Y}^\tau :=\overline{Y}^\tau
\cap \h$ without further explanation.
\medskip

\noindent{\bf Simple homeomorphisms.} In the paper \cite{Gl-W 2}, a new
interesting class of homeomorphisms was defined. By definition, $S\in \h$
is {\it simple} if it satisfies the following conditions.

(i) There exist non-empty clopen subsets $F_j$ and integers $r_j\ge 1,\
j=1,...,k,$ such that the collection  $\{S^iF_j : i=0,1,...,r_j,\
j=1,...,k\}$ is pairwise disjoint and $S$ has period $r_j$ on $F_j$.

(ii) There exist clopen subsets $C_s,\ s=1,...,l,$ and, for each $s$, two
disjoint periodic orbits $(y_s^+, Sy_s^+,...,S^{q^+_s-1}y_s^+)$ and
$(y_s^-, Sy_s^-,...,S^{q^-_s-1}y_s^-)$ such that the sets $(S^nC_s : n\in
\Z, s=1,...,l)$ are pairwise disjoint and spiral towards the periodic
orbits of $y_s^+$ and  $y_s^-$, that is $\lim_{n\to \pm \infty}{\rm
dist}(S^nC_s, S^ny_s^{\pm}) = 0$

(iii) The space $\Om$ may be represented as
\begin{equation}\label{simple}
\Om = \bigcup_{j=1}^k\bigcup_{i=0}^{r_j-1} S^iF_j\cup \
\bigcup_{s=1}^l\bigcup_{n\in \Z} S^nC_s\ \cup\ \bigcup_{s=1}^l
\left[(y_s^+,...,S^{q^+_s-1}y_s^+) \cup (y_s^-,
...,S^{q^-_s-1}y_s^-)\right]
 \end{equation}

It follows from the definition that every simple homeomorphism $S$ contains
a periodic part and there may exists a clopen set $C$ such that $TC
\subsetneqq C$ or $C \subsetneqq TC$.

In \cite{Gl-W 2}, the authors proved the following important result which
was used to establish the Rokhlin property in $\h$.

\begin{theorem}\label{SimpleIsDense} The set of simple homeomorphisms is
dense in $(\h,p)$.
\end{theorem}

%
%  PERIODIC AND APERIODIC HOMEOMORPHISMS
%

\noindent{\bf Periodic and aperiodic homeomorphisms.} As an immediate
consequence of the Rokhlin lemma, we can find some topological properties
of the sets of periodic and aperiodic homeomorphisms. The following two
theorems was proved in \cite{B-D-K 2}.

\begin{theorem}[aperiodic]\label{Aperiodic_Homeo}
$(1)$ $\ap$ is closed and nowhere dense in $(\h,\tau)$.\\
$(2)$ $\ap$ is dense in $\h$ with respect to the topology $p \ (\sim D)$.
\end{theorem}

To prove (1), observe that  $\ap$ is closed and $\h\setminus \ap$ is dense
by Theorem \ref{RokhlinLemmaCantor}(3). The second assertion follows from
the fact that every simple homeomorphism can be approximated in the
topology $p$ by aperiodic ones. Then, we use Theorem \ref{SimpleIsDense} to
complete the proof.

\begin{theorem}[periodic]\label{PeriodicHomeo}
$(1)$ $\overline{\per_0}^\tau = \h$.\\
$(2)$ The set $\overline{\per}^p$ is a proper subset in $\h$ containing
minimal homeomorphisms.\\
$(3)$ $\overline{\per}^p=\overline{\per_0}^p$.
\end{theorem}

Since the first result is new we give the proof of it. Observe that the
fact that $\overline{\per_0}^\tau\supset \ap$ follows from the Rokhlin
lemma (Theorem \ref{RokhlinLemmaCantor}). Now, we give an idea how to
approximate any $T\in\h$ by elements from $\per_0$. For $T\in \h$, consider
a neighborhood $U = U(T;\mu_1,\ldots,\mu_k;\e)$. Take the canonical
partition of $\Om=X_\infty\cup\bigcup_{n\geq 1}X_n$ where $X_n$ is the set
of points with period $n$ and $X_\infty$ consists of all aperiodic points.
Notice that every $X_n=\bigcup_{i=0}^{n-1} T^i X_n^0,\ 0 < n< \infty,$ is a
$T$-tower. By Proposition \ref{K-R partitions}, we can find a K-R partition
$\Xi=\{\xi_1,\xi_2,\ldots\}$ of $X_\infty$ such that
$\mu_l(T^{-1}B(\Xi))<\e/4$ for $l=1,\ldots,k$ where $B(\Xi)$ is the base of
$\Xi$. Find $n$ and $m$ such that $\mu_l(\xi_1\cup\ldots\xi_n\cup
X_1\cup\ldots\cup X_m)>1-\e/4$. Then, by regularity of the measures
$\mu_l$, we can find in each $T$-tower $\xi_i$ and $X_j$, a closed subtower
$\xi'_i$ and $X'_j$, resp., such that $\mu_l(\xi_i\setminus
\xi'_i)<\e(4n)^{-1}$ and $\mu_l(X_j\setminus X'_j)<\e(4m)^{-1},\ i=1,...,n;
j=1,...,m$. Thus, we have obtained a finite number of disjoint closed
$T$-towers such that $\mu_l(\xi'_1\cup\ldots\cup\xi'_n\cup X'_1\cup\ldots
X'_m)>1-\e/2$ and $\mu_l(T^{-1}B(\xi'_1\cup\ldots\cup\xi'_n))<\e/2$. For
every $\xi'_i $ and $X'_j$, we find a clopen  $T$-tower $\xi_j^{co}$ and
$X_j^{co}$, resp.,  such that $\xi'_i$ and $X'_j$ are subtowers of
$\xi_j^{co}$ and $X_j^{co}$ and $\mu_l(\bigcup_{i=1}^n
T^{h(\xi_i^{co})}B(\xi_i^{co})) < \e/2$. Moreover, we can choose all these
clopen $T$-towers pairwise disjoint. Let $Q$ be the union of these clopen
$T$-towers. Define $P\in\h$ as follows: if $x\notin Q$, set $Px=x$; if
$x\in Q$ and does not lie in the top level of a tower, set $Px=Tx$, and for
$x$ lying in the top level of a $T$-tower from $Q$ (say $\xi$), set
$Px=T^{-h(\xi)+1}x$, where $h(\xi)$ is the height of $\xi$. Observe that
$P\in U\cap\per_0$.

Statement (2) follows from the following observation: every dissipative
homeomorphism does not belong to $\overline{\per}^p$. Recall that a
homeomorphism $T$ is  {\it dissipative} if there is a clopen set $F$ such
that either $TF\subsetneqq F$ or $F\subsetneqq TF$. The fact that $\mathcal
Min \subset \overline{\per_0}^p$ was proved in \cite{B-D-K 2} where the
closure of $\per_0$ in $p$ was found.
\\

%
% %%%%%%%%%%%%%%%%%   ODOMETERS, MINIMAL AND MIXING
%

\noindent{\bf Odometers, minimal, and mixing homeomorphisms.} This section
contains our results on the description of closures of some natural subsets
in $\h$ with respect to the topologies $\tau$ and $p$. The notation
$\mathcal Od$, $\mathcal Min$, $\mathcal Mix$ was introduced in Section 1.
Denote also by $\mathcal Tt$ the set of topologically transitive
homeomorphisms. Clearly, $\mathcal Od\subset \mathcal Min \subset \mathcal
Tt\cap \ap \subset\ap$. As in ergodic theory, one can define the set
$\mathcal R(1)$ of homeomorphisms of rank 1 but it is not hard to show that
this set coincides with $\mathcal Od$. The following definition was
introduced in \cite{B-D-K 2}: A homeomorphism $T$ is called {\it moving} if
for every proper clopen set $F$ each of the sets $F\setminus TF$ and
$TF\setminus F$ is not empty. We denote the class of moving homeomorphisms
by $\mathcal Mov$. It is clear that $\per\cap\mathcal Mov=\emptyset$ and
$\mathcal Mov \supset \mathcal Tt \supset Min$.

The next results were proved in \cite{B-D-K 2}.

\begin{theorem}\label{ClassesHomeom}
$(1)$ The sets $\overline{ \mathcal Min}^{p}$ and $\overline{\mathcal
Mix}^p$ are nowhere dense in $(\h,p)$.\\
$(2)$ $\mathcal Od$ is a dense $G_\delta$-set in
$\overline{\mathcal Min}^p$.\\
$(3)$ $\mathcal Od$ is an  $F_{\sigma\delta}$-set in $\overline{\mathcal
Min}^\tau$.\\
$(4)$ $\overline{\mathcal Od}^\tau=  \overline{\mathcal Min}^\tau=
\overline{\mathcal Tt}^\tau \cap \ap = \ap$ and $\overline{\mathcal
Od}^p=\overline{\mathcal Min}^p= \overline{\mathcal Tt}^p= \mathcal Mov$.\\
$(5)$ $\mathcal Mix\cap  \overline{\mathcal Min}^p$ is nowhere dense  in
$(\overline{\mathcal Min}^p,p)$.
\end{theorem}

In \cite{Gl-W 2}, the following interesting theorem was proved.

\begin{theorem}\label{topEntropy} $(1)$ For a compact metric space $K$,
the set of homeomorphisms having zero entropy is a $G_\delta$ subset in
$Homeo(K)$ in the topology of uniform convergence.\\
$(2)$ The set of homeomorphisms having zero entropy is a dense $G_\delta$
subset in $(\h, p)$.\\
$(3)$ The set of homeomorphisms of the Hilbert cube having infinite
topological entropy contains a dense $G_\delta$ subset in the topology of
uniform convergence. This is also true for a finite dimensional cube.
\end{theorem}

%
%%%%%%%%%%%%%%%%%%%%%%%%% FULL GROUPS
%

\noindent{\bf Subgroups of $\h$ generated by a homeomorphism.} In this
section, we discuss topological properties of full groups, normalizers, and
centralizers of Cantor systems.

Let $T$ be an aperiodic homeomorphism. It is easy to prove that the full
group $[T]$ is closed and nowhere dense in $(\h,\tau)$. On the other hand,
$[T]$ is not closed in $(\h,p)$, but $[T]$ is still nowhere dense.

The following theorem asserts that the full group of a minimal
homeomorphism $T$ is separable since the topological full group is a dense
countable subgroup in $[T]$ in the topology $\tau$. Notice that this result
is not true, in general, with respect to the topology $p$. We find the
class of minimal homeomorphisms for which the $p$-closures of $[[T]]$ and
$[T]$ coincide. We say that a minimal homeomorphism $T$ is {\it saturated}
if for any two clopen subsets $A$ and $B$ with $\mu(A)=\mu(B)$ for every
$T$-invariant measure $\mu$ there is $\gamma\in[[T]]$ such that
$\gamma(A)=B$ \cite{B-K 1}. Observe that the class of saturated
homeomorphisms contains odometers and does not contain the Chacon
homeomorphisms. By Theorem \ref{ClassesHomeom}(2), the class of saturated
homeomorphisms is generic in $\mathcal Mov$ with respect to $p$.

\begin{theorem}\label{FullGroupsCantor}
Let $T$ be a minimal homeomorphism of $\Om$.\\
$(1)$ The topological full group $[[T]]$ is $\tau$-dense in $[T]$.\\
$(2)$ $\overline{[T]}^p=\{S\in \h : \mu\circ S=\mu \mbox{ for all }\mu\in
M_1(T)\}$.\\
$(3)$ $T$ is saturated if and only if $\overline{[[T]]}^p=\overline
{[T]}^p$.\\
$(4)$ If $T$ is of rank 1 (odometer), then $C(T)= \overline{\{T^n : n \in
\Z\}}^p$.
\end{theorem}

The first and third statements are proved in \cite{B-K 1}. The second
result is taken from \cite{G-P-S 2}, and (4) is a rather simple observation
noticed in \cite{B-D-K 2}.
\smallskip

For a homeomorphism  $T\in\h$, let $\per_0(T)$ denote the set $\per_0\cap
[[T]]$. Although $\per_0$ is not dense in $(\h,p)$,  the following result
holds. The first statement follows from Theorem \ref{PeriodicHomeo}(1), and
(2) was proved in \cite{B-D-K 2}.

\begin{theorem} $(1)$ Let $T$ be a homeomorphism  of $\Om$, then
 $\overline{\per_0(T)}^\tau=[T]$.\\
$(2)$ Let $T$ be a minimal homeomorphism of $\Om$, then
$\overline{\per_0(T)}^p\supset [[T]].$
\end{theorem}

Let $T$ be a homeomorphism of a Cantor set $\Om$. Then $T$ can be also
considered as a Borel automorphism of $(\Om, \B)$. Therefore, one can
define two full groups $[T]_C$ and $[T]_B$:
$$
[T]_C = \{S\in \h : Sx \in \{T^nx\ |\ n\in \Z\} \ \forall x\in \Om\},
$$
$$
[T]_B = \{S\in Aut(\Om,\B): Sx \in \{T^nx\ |\ n\in \Z\} \ \forall x\in
\Om\}
$$
where $C$ and $B$ stand for Cantor and Borel.  Obviously, $[[T]]_C\subset
[T]_C \subset [T]_B$ and $[T]_B$ is closed in $Aut(\Om,\B)$ with respect to
$\tau$. Recall that by Theorem \ref{FullGroupsCantor} (1) we have that for
a minimal homeomorphism $T$,
\begin{equation}
\overline{[[T]]_C}^\tau \cap \h = [T]_C.
\end{equation}

The problem of finding of the closure of $[[T]]_C$ in $Aut(\Om, \B)$ with
respect to $\tau$ was solved in \cite{B-D-K 2} where the following result
was obtained.

\begin{theorem} Let $T$ be a minimal homeomorphism of a Cantor
set $\Om$. Then $\overline{[[T]]_C}^\tau = [T]_B$.
\end{theorem}

This theorem has been recently generalized to any aperiodic homeomorphism
$T$ by Medynets.

We finish this section by considering a topology on the normalizer $N[T]$
of $T\in \h$. By definition, the (topological) normalizer of a
homeomorphism $T\in \h$ is the group $N[T]=\{S\in \h : S(Orb_T(x))=
Orb_T(Sx)\mbox{ for all }x\in X\}$. In \cite{B-K 2}, the following topology
on $N[T]$ was introduced: for $R\in N[T]$, Borel probability measures
$\mu_1,\ldots,\mu_k$, and $\e>0$ define $V(R;\mu_1,\ldots,\mu_k;\e) =\{P\in
N[T] : D(R, P)<\e \ \mbox{ and }\mu_i(E(RT^nR^{-1}, PT^nP^{-1}))<\e,\;
i=1,\ldots,k, n\in \Z \}$. The sets $\{V(R;\mu_1,\ldots,\mu_k;\e)\}$ form a
base of a topology $\lambda$. The following result is due to \cite{B-K 2}.

\begin{theorem} $(N[T],\lambda)$ is a complete separable topological group.
\end{theorem}

%%%%%%%%%%%%%%%%%%%%%Comparison%%%%%%%%%%%%

\section{Comparison of topological properties of transformation groups.}

In this section, we  discuss the results from Sections 2 - 4. The reader
can see that many of the formulated theorems are similar for measurable,
Borel, and Cantor dynamics. This observation shows that the three dynamics
have common properties in the topologies we have studied.

We first consider the groups $\Aut$, $\aut$, and $\h$ with respect to the
uniform topology and summarize their topological properties in the
following table.

%
%%%%%%%%%%%%%%%%%%%%  THE TABLE FOR UNIFORM TOPOLOGY
%
\begin{center}
{\bf Uniform topology}
\end{center}
\noindent\begin{tabular}
{|p{2.2cm}|p{1.9cm}|p{1.6cm}|p{1.6cm}|p{1.6cm}|p{1.9cm}|}\hline {Group} &
$\Aut$ & \multicolumn{2}{c|}{$\aut$}

& $\widehat{Aut}\bs$ & $\h$\\ \hline

{Topology} & $d_u$ & $\tau$ &  $\tau_0$ & $\widehat\tau$ & $\tau$\\
\hline

Connectedness & Yes & totally disconnected & Yes & Yes & totally dicsonected\\
\hline

Path-connectedness & Yes & No & Yes & Yes & No\\ \hline

Contractibility & Unknown & No & Unknown & Unknown & No\\ \hline

{Rokhlin lemma and Rokhlin property}   &    holds &

 holds & holds & holds & holds\\
\hline

$\per$ & dense &  dense & dense & dense & dense\\ \hline

$\ap$  & closed nowhere dense & closed nowhere dense & nonclosed
nowhere dense & closed nowhere dense & closed nowhere dense\\
\hline

Dense subsets & $\overline{\mathcal Erg}^{d_u}=\ap$ $\overline{\mathcal
Od}^{d_u}=\ap$ & $\overline{\mathcal Od}^\tau=\ap$ & $\overline{\mathcal
Od}^\tau=\ap$ & $\overline{\mathcal Od}^\tau=\ap$ &
$\overline{\mathcal Min}^\tau=\ap$ $\overline{\mathcal Od}^\tau=\ap$\\
\hline
\end{tabular}
\medskip

It is seen that if the topological groups from the table are connected then
they are path-connected. The Rokhlin lemma  and Rokhlin property hold for
all these groups. Recall that by the Rokhlin property in the uniform
topology we mean the density of conjugates of any aperiodic transformation
in $\ap$. Another common property for these groups is the density of
periodic transformations. This fact allows one to approximate every
transformation by periodic ones. On the other hand, the set of aperiodic is
closed and nowhere dense in all the dynamics. However, inside of the set of
aperiodic transformations, we can find dense subsets consisting of
transformations which are of \textquotedblleft ergodic" type, i.e.
odometers, minimal homeomorphisms. As seen from the table and obtained
results,  the best analogue of the space $(\Aut, d_u)$ in Borel dynamics is
apparently the quotient group $(\widehat{Aut}(X,\B), \widehat{\tau})$ (or
$(\aut, \tau_0))$. Note that we have not considered properties of $\h$ with
respect to $\tau_0$. We think it may be interesting to do this.

\begin{center}
\textbf{Weak topology}
\end{center}

A similar comparison of these groups  with respect to the weak topologies,
$d_w$ and $p$, shows that they do not have  much in common in contrast to
the uniform topology. Nevertheless, we can point out some common features.
Firstly, the set of aperiodic transformations is dense for each of these
groups in the weak topology. Secondly, the Rokhlin property holds in
$\Aut$, $\widehat{Aut}(X,\B)$, and $\h$. Moreover, the following two facts
hold: (i) the  transformations of zero entropy are typical for $\Aut$ and
$\h$, (ii) rank 1 automorphisms are typical in $Aut_0(X,\B,\mu)$ and
odometers (or rank 1 homeomorphisms) are typical in the space $\mathcal
Mov$, the closure of all minimal homeomorphisms. Based on Theorem
\ref{ClassesHomeom}, one can suppose that the set $\mathcal Mov$ endowed
with the topology $p$ is an analogue of the space $(\Aut, d_w)$. On the
other hand, the group $\Aut$ is contractible in the weak topology, but
$(\widehat{Aut}\bs,\widehat p)$ and $(\h,p)$ are totally disconnected
topological spaces.

For a Cantor set $\Om$, the topology of uniform convergence defined by the
metric $D$ and the topology $p$ are equivalent. If $(X, d)$ is a compact
metric space then the topology $D$ and $p$ on $\aut$ are different. It is
interesting to compare topological properties of $(\h,p)$ and $(\aut, D)$.
We see that the set $\ap$ of periodic transformations is dense in $(\Aut,
d_w)$, $(\h, p)$, and ``almost'' dense in $(\aut, D)$. We hope that the
further study of possible relations between these topologies will lead to
new interesting results.

%%%%%%%%%%%%%%%%%%%%%%%%%%%%%%%%%%%%%%%

\bibliographystyle{amsalpha}

\end{document}